\documentclass[a4paper,11pt]{article}
\usepackage[top=2.5cm,bottom=2.5cm,left=2.2cm,right=2.2cm]{geometry}
\usepackage[T1]{fontenc}
\usepackage{authblk}
\usepackage{lmodern}
\usepackage{amsmath,amssymb,amsfonts,amsthm,mathtools}
\usepackage{microtype}
\usepackage{indentfirst}
\usepackage{enumitem}
\usepackage{hyperref}
\hypersetup{colorlinks=true,linkcolor=blue,citecolor=blue,urlcolor=blue}
\usepackage{bookmark}

\allowdisplaybreaks
\makeatletter
\renewcommand{\@seccntformat}[1]{{\csname the#1\endcsname}{\normalsize .}\hspace{.5em}}
\makeatother

\numberwithin{equation}{section}

\setlist{nosep}

\theoremstyle{plain}
\newtheorem{thm}{Theorem}[section]
\newtheorem{lem}[thm]{Lemma}

\newtheorem{pb}[thm]{Problem}
\newtheorem{claim}{Claim}
\theoremstyle{definition}
\newtheorem{conj}[thm]{Conjecture}
\theoremstyle{remark}

\def\b{\color{blue}}
\newcommand{\Cdir}{\overrightarrow{C}}
\newcommand{\ex}{\operatorname{ex}}

\newcommand{\keywords}[1]{\par\vspace{6pt}\noindent\textbf{Keywords: }#1\par}
\newcommand{\MSC}[1]{\par\vspace{3pt}\noindent\textbf{MSC(2020): }#1\par}

\hypersetup{
    colorlinks=true,
    linkcolor=blue,
    citecolor=red,
    urlcolor=magenta,
}

\begin{document}
\baselineskip=0.23in

\title{\bf Directed triangles in oriented graphs without a directed cycle of length $3t$}
\author{Shuchao Li}
\author{Mingli Wang\footnote{Corresponding author.\\
\hspace*{2em}E-mail: wmlmath@163.com}}
\affil{School of Mathematics and Statistics, and Hubei Key Lab--Math. Sci.,\linebreak Central China Normal University, Wuhan 430079, China}
\date{\today}

\maketitle

\begin{abstract}
For integers $k,\ell\ge 3$, let $\operatorname{ex}(n,\Cdir_k,\Cdir_\ell)$ denote the maximum number of directed cycles of length $k$ in an $n$-vertex oriented graph containing no directed cycles of length $\ell$. Grzesik, Jaworska, Kielak, Kuc and \'Slusarczyk~\cite{GJ} systematically studied the generalized Tur\'an problem for directed cycles and conjectured that, for every fixed integer $t\ge 2$, $\operatorname{ex}(n,\Cdir_3,\Cdir_{3t})=\frac{t-1}{4}n^2+o(n^2)$. 
They proved this asymptotic formula for $t=2$. In this paper, we prove it for every fixed integer $t\ge 3$, thereby confirming the conjecture.
\end{abstract}
\MSC{05C35; 05C20}

\section{\normalsize Introduction}

One of the central problems in extremal graph theory is the classical Tur\'an problem, which asks for the maximum number of edges in an $n$-vertex graph containing no copy of a fixed graph $F$. 
This maximum, denoted by $\operatorname{ex}(n,F)$, is called the \textit{Tur\'an number} of $F$. The study of this problem originated with Mantel's theorem~\cite{Mantel} on triangle-free graphs, which was later generalized by Tur\'an~\cite{Turan} to forbidden complete graphs. The Erd\H{o}s--Stone theorem~\cite{ErdosStone} determines the asymptotic value of $\operatorname{ex}(n,F)$ for every fixed non-bipartite graph $F$. In contrast, determining Tur\'an numbers for bipartite graphs is generally more difficult. Bondy and Simonovits~\cite{BondySimonovits} established a fundamental upper bound on the Tur\'an numbers of even cycles.

A natural extension of the classical Tur\'an problem replaces the number of edges by the number of copies of a fixed graph. For graphs $H$ and $F$, the \textit{generalized Tur\'an number} $\operatorname{ex}(n,H,F)$ is the maximum number of copies of $H$ in an $n$-vertex $F$-free graph. This framework was systematically developed by Alon and Shikhelman~\cite{AlonShikhelman2016}, whose results cover several fundamental classes of counted and forbidden graphs. One important case is when both $H$ and $F$ are cycles. The Erd\H{o}s pentagon conjecture, proved independently by Grzesik~\cite{GrzesikPentagons} and by Hatami, Hladk\'y, Kr\'a\v{l}, Norine and Razborov~\cite{HatamiPentagons}, states that every $n$-vertex triangle-free graph contains at most $(n/5)^5$ copies of $C_5$. Gerbner, Gy\H{o}ri, Methuku and Vizer~\cite{GGMV} studied generalized Tur\'an problems for even cycles, while Beke and Janzer~\cite{BekeJanzer2024} investigated the corresponding problems for odd cycles.

Directed analogues of Tur\'an-type problems have a separate history. Brown, Erd\H{o}s and Simonovits~\cite{BrownErdosSimonovits1973} developed an asymptotic structural theory for maximizing the number of arcs in digraphs subject to forbidden-subdigraph conditions. We distinguish between digraphs (without loops or multiple arcs), which may contain both arcs $uv$ and $vu$, and oriented graphs, which contain at most one of them.
Zhou and Li~\cite{ZhouLi2023} obtained exact results for Tur\'an numbers of directed paths and directed cycles in digraphs.

For Tur\'an problems on oriented graphs, the distinction becomes especially significant when the forbidden graph contains a directed cycle. 
A transitive tournament has $\binom{n}{2}$ arcs and contains no directed cycle, so the corresponding arc extremal number is already $\binom{n}{2}$. This observation explains the emphasis on acyclic forbidden graphs in classical oriented Tur\'an theory. 
In that setting, Grzesik, Jaworska, Kielak, Novik and \'Slusarczyk~\cite{GrzesikEtAl2024} studied compressibility and its role in determining asymptotic extremal arc densities.
Counting directed subgraphs provides a natural way to formulate nontrivial extremal questions when the forbidden graph is a directed cycle.

For oriented graphs $H$ and $F$, let $\operatorname{ex}(n,H,F)$ denote the maximum number of copies of $H$ in an $n$-vertex $F$-free oriented graph. Let $\Cdir_\ell$ denote the directed cycle of length $\ell$.
Grzesik, Jaworska, Kielak, Kuc and \'Slusarczyk~\cite{GJ} systematically studied generalized Tur\'an problems for directed cycles in oriented graphs. They determined the order of magnitude of $\operatorname{ex}(n,\Cdir_k,\Cdir_\ell)$ for all fixed distinct integers $k,\ell\ge3$, proving that $\operatorname{ex}(n,\Cdir_k,\Cdir_\ell)=\Theta(n^k)$ when $k\nmid\ell$, and $\operatorname{ex}(n,\Cdir_k,\Cdir_\ell)=\Theta(n^{k-1})$ when $k\mid\ell$. They also obtained asymptotically sharp results in several cases where $k\nmid\ell$. In particular, for directed triangles, they proved that
$\ex(n,\Cdir_3,\Cdir_\ell)=\frac{n^3}{27}+O(n^2)$ for every fixed integer $\ell\ge4$ not divisible by three~\cite[Theorem~21]{GJ}. However, when $k\mid\ell$, determining the leading term remains open in general. For directed triangles, they 
proposed the following conjecture.

\begin{conj}[\cite{GJ}, Conjecture 35]\label{conj:main}
For every fixed integer $t\ge2$, 
\begin{align}\label{e1.1}
\ex(n,\Cdir_3,\Cdir_{3t})=&\frac{t-1}{4}n^2+o(n^2).
\end{align}
\end{conj}
Grzesik, Jaworska, Kielak, Novik and \'Slusarczyk~\cite{GJ} showed that \eqref{e1.1} holds for $t=2$. Our main result establishes \eqref{e1.1} for every fixed integer $t\ge3$, thereby resolving Conjecture~\ref{conj:main} completely.

\begin{thm}\label{thm:main}
For every fixed integer $t\ge3$, $\ex(n,\Cdir_3,\Cdir_{3t})=\frac{t-1}{4}n^2+o(n^2)$.
\end{thm}

\noindent{\bf Organization.}
Section~\ref{s2} introduces some notation and collects preliminary results needed for the proof. Section~\ref{s3} contains the proof of Theorem~\ref{thm:main}. Finally, Section~\ref{s4} presents concluding remarks and discusses related questions concerning longer directed cycles.

\section{\normalsize Some notation and preliminaries}\label{s2}

For an oriented graph $G$, let $V(G)$ and $A(G)$ denote its vertex set and arc set, respectively. We write $n(G)=|V(G)|$ and $e(G)=|A(G)|$. 
An arc from $u$ to $v$ is denoted by $uv$. For each vertex $v\in V(G)$, define $N_G^+(v)=\{u\in V(G):vu\in A(G)\}$ and $N_G^-(v)=\{u\in V(G):uv\in A(G)\}$. The out-degree and in-degree of $v$ are $d_G^+(v)=|N_G^+(v)|$ and $d_G^-(v)=|N_G^-(v)|$, respectively. For $S\subseteq V(G)$, write $N_S^\pm(v)=N_G^\pm(v)\cap S$. We omit the subscript $G$ when there is no ambiguity. For $X\subseteq V(G)$, let $G[X]$ denote the subgraph of $G$ induced by $X$, and write $G-X=G[V(G)\setminus X]$. A matching in an oriented graph is a set of pairwise vertex-disjoint arcs.

For an undirected graph $G$, let $E(G)$ denote its edge set, and write $n(G)=|V(G)|$ and $e(G)=|E(G)|$. An edge between $u$ and $v$ is denoted by $uv$. For $v\in V(G)$, let $N_G(v)$ and $d_G(v)$ denote its neighbourhood and degree, respectively. We write $\Delta(G)$ for the maximum degree of $G$. The undirected cycle of length $\ell$ is denoted by $C_\ell$. For graphs $H$ and $G$, either both undirected or both oriented, let $N(H,G)$ denote the number of subgraphs of $G$ isomorphic to $H$; copies are not required to be induced.

The asymptotic notation $O(\cdot)$, $o(\cdot)$, and $\Theta(\cdot)$ is used as $n\to\infty$. The subscript in $O_t(\cdot)$ indicates that the implied constant may depend on $t$.

We recall the following extremal and counting results, which will be used in the proof.

\begin{lem}[\cite{ErdosStone}]\label{lem:erdos-stone}
For every fixed graph $F$ with $\chi(F)\ge2$,
\[
\operatorname{ex}(n,F)
=\left(1-\frac{1}{\chi(F)-1}+o(1)\right)\frac{n^2}{2}.
\]
\end{lem}

\begin{lem}[\cite{BondySimonovits}]\label{lem:BS}
For every integer $t\ge2$, $\operatorname{ex}(n,C_{2t})\le 100t\,n^{1+1/t}.$
\end{lem}

\begin{lem}[\cite{GGMV}]\label{lem:c4c2t}
For every fixed integer $t\ge2$,
\[
\operatorname{ex}(n,C_4,C_{2t}) =(1+o(1))\frac{(t-1)(t-2)}{4}n^2 = O_t(n^2).
\]
\end{lem}

\begin{lem}[\cite{Fox}]\label{lem:removal}
For every fixed undirected graph $F$ and every $\varepsilon>0$, there exists $\delta>0$ such that every $n$-vertex undirected graph containing at most $\delta n^{n(F)}$ copies of $F$ can be made $F$-free by deleting at most $\varepsilon n^2$ edges.
\end{lem}

For undirected graphs $H$ and $G$, a \textit{homomorphism} from $H$ to $G$ is a map $\varphi:V(H)\to V(G)$ such that $\varphi(u)\varphi(v)\in E(G)$ whenever $uv\in E(H)$. We write $\hom(H,G)$ for the number of such homomorphisms. In particular, a homomorphic $2k$-cycle in an undirected graph $G$ is a sequence $(x_1,\ldots,x_{2k})$ of vertices of $G$ such that $x_ix_{i+1}\in E(G)$ for every $i\in[2k]$, where indices are taken modulo $2k$. The vertices $x_1,\ldots,x_{2k}$ need not be distinct.

\begin{lem}[\cite{Janzer}]\label{lem:janzer}
Let $k\ge2$ be an integer, and let $G$ be an undirected graph on $n$ vertices. Let $\sim$ be a symmetric binary relation defined over $V(G)$ such that, for every $u,v\in V(G)$, $|\{w\in N_G(v):u\sim w\}|\le s.$
Then the number of homomorphic $2k$-cycles $(x_1,\ldots,x_{2k})$ in $G$ such that $x_i\sim x_j$ for some $i\ne j$ is at most
$$
32k^{3/2}s^{1/2}\Delta(G)^{1/2}n^{1/(2k)}\hom(C_{2k},G)^{1-1/(2k)}.
$$
\end{lem}

\begin{lem}[\cite{HornJohnson2013}]\label{lem:Rayleigh}
Let $G$ be an undirected graph on $n$ vertices with adjacency matrix $A(G)$, and let $\rho(G)$ denote the spectral radius of $A(G)$, that is, the maximum absolute value of its eigenvalues. Then, for every nonzero vector $\mathbf{x}\in\mathbb{R}^n$,
\[
\rho(G)\geq \frac{\mathbf{x}^{\mathsf T}A(G)\mathbf{x}} {\mathbf{x}^{\mathsf T}\mathbf{x}}.
\]
Equality holds if and only if $\mathbf{x}$ is an eigenvector of $A(G)$ corresponding to $\rho(G)$.
\end{lem}

\section{\normalsize Proof of Theorem~\ref{thm:main}}\label{s3}

To prove Theorem~\ref{thm:main}, we first establish three auxiliary lemmas.

\begin{lem}\label{lem:dirremoval}
For every $\eta>0$, there exists $\xi>0$ such that every $n$-vertex oriented graph containing at most $\xi n^3$ copies of $\Cdir_3$ can be made $\Cdir_3$-free by deleting at most $\eta n^2$ arcs.
\end{lem}

\begin{proof}
Fix $\eta>0$, and let $\delta>0$ be the constant given by Lemma~\ref{lem:removal} with $F=C_3$ and $\varepsilon=\eta/9$. Set $\xi=9\delta$.
Let $G$ be an $n$-vertex oriented graph containing at most $\xi n^3$ copies of $\Cdir_3$. Construct an undirected tripartite graph $B$ with pairwise disjoint vertex classes $V_i=\{v_i:v\in V(G)\}$ for $i=1,2,3$. For each arc $xy\in A(G)$, add the three edges $x_1y_2$, $x_2y_3$, and $x_3y_1$ to $B$, and no others.

Since $B$ is tripartite, every triangle in $B$ has exactly one vertex in each of $V_1,V_2,V_3$, say $x_1,y_2,z_3$. By construction, these vertices form a triangle in $B$ if and only if $xy,yz,zx\in A(G)$. Moreover, each copy of $\Cdir_3$ in $G$ gives rise to exactly three triangles in $B$.
Hence $B$ contains at most $3\xi n^3=\delta(3n)^3$ triangles.
Since $B$ has $3n$ vertices, Lemma~\ref{lem:removal} implies that there exists a set $R\subseteq E(B)$ with
$|R|\le \frac{\eta}{9}(3n)^2=\eta n^2$ such that $B-R$ is triangle-free.
For each edge in $R$, delete the corresponding arc from $G$. Since each edge of $B$ corresponds to a unique arc of $G$, at most $|R|\le\eta n^2$ arcs are deleted from $G$.

If the resulting oriented graph contained a directed triangle, its three corresponding edges would form a triangle in $B-R$, a contradiction. Thus the resulting oriented graph is $\Cdir_3$-free.
\end{proof}

\begin{lem}\label{lem:dsq}
Let $t\ge3$ be a fixed integer. Then the following statements hold.
\begin{enumerate}
\item[\textup{(i)}] There exists a constant $c_t>0$ such that every $n$-vertex $C_{2t}$-free bipartite graph $\Gamma$ satisfies $\sum_{v\in V(\Gamma)}d_\Gamma(v)^2\le c_tn^2$.
\item[\textup{(ii)}] For every $\varepsilon>0$, there exists $\delta=\delta(t,\varepsilon)>0$ such that every $n$-vertex $C_{2t}$-free bipartite graph $\Gamma$ with $\Delta(\Gamma)\le\delta n$ satisfies $\sum_{v\in V(\Gamma)}d_\Gamma(v)^2\le\varepsilon n^2$.
\end{enumerate}
\end{lem}
\begin{proof}
Choose $c_t\ge5$ sufficiently large, depending only on $t$, such that $N(C_4,\Gamma)\le c_tn^2/16$ for every $n$-vertex $C_{2t}$-free bipartite graph $\Gamma$, and $\operatorname{ex}(m,C_t)\le 3m^2/8$ for every integer $m\ge c_t$. Such a choice is possible by Lemma~\ref{lem:c4c2t}, together with Lemma~\ref{lem:erdos-stone} when $t$ is odd and Lemma~\ref{lem:BS} when $t$ is even. We fix this constant throughout the proof.

We first establish (i). Let $(A,B)$ be the bipartition of $\Gamma$, and write $d(v)=d_\Gamma(v)$. For distinct vertices $x,y\in A$, let $q(x,y)=|N_\Gamma(x)\cap N_\Gamma(y)|$. Counting copies of $C_4$ by their two vertices in $A$, we have
\[
\sum_{\{x,y\}\subseteq A}\binom{q(x,y)}2=N(C_4,\Gamma)\le\frac{c_t}{16}n^2.
\]
Since $q\le1+2\binom q2$ and $q^2=q+2\binom q2$ for every nonnegative integer $q$, it follows that
\begin{equation}\label{eq:codegree}
\sum_{\{x,y\}\subseteq A}q(x,y)\le\left(\frac12+\frac{c_t}{8}\right)n^2,\qquad \sum_{\{x,y\}\subseteq A}q(x,y)^2\le c_tn^2.
\end{equation}
On the other hand, counting pairs of vertices in $A$ with a common neighbour in $B$ gives
\[
\sum_{b\in B}\binom{d(b)}2=\sum_{\{x,y\}\subseteq A}q(x,y).
\]
Consequently,
\[
\sum_{b\in B}d(b)^2=e(\Gamma)+2\sum_{\{x,y\}\subseteq A}q(x,y)\le\left(\frac54+\frac{c_t}{4}\right)n^2.
\]
The same bound holds with $A$ and $B$ interchanged. Hence
\[
\sum_{v\in V(\Gamma)}d(v)^2\le\left(\frac52+\frac{c_t}{2}\right)n^2\le c_tn^2.
\]

We now prove (ii). Fix $\varepsilon>0$. Choose $\eta>0$ such that $0<\eta\le\min\left\{1,\frac{\varepsilon^2}{1024c_t}\right\}$.
By Lemma~\ref{lem:removal} applied to $C_t$ with parameter $\eta$, choose $\rho>0$ such that every graph on $n$ vertices with at most $\rho n^t$ copies of $C_t$ can be made $C_t$-free by deleting at most $\eta n^2$ edges. Choose $\delta>0$ such that
\[
0<\delta\le\min\left\{\frac12,\frac{\rho}{t(t-1)c_t},\frac{\varepsilon}{2c_t}\right\}.
\]

Let $\Gamma$ be an $n$-vertex $C_{2t}$-free bipartite graph with bipartition $(A,B)$ and $\Delta(\Gamma)\le\delta n$. Define a graph $F$ on $A$ by joining two distinct vertices whenever they have a common neighbour in $\Gamma$. For each edge $xy\in E(F)$, choose a common neighbour $b_{xy}\in B$.

Consider a cycle $x_1x_2\cdots x_tx_1$ in $F$. If the chosen common neighbours of its edges were all distinct, they would form a copy of $C_{2t}$ in $\Gamma$ together with $x_1,\ldots,x_t$, a contradiction. Thus two distinct edges of the cycle have the same chosen common neighbour $b$. If these edges are adjacent, three vertices of the cycle lie in $N_\Gamma(b)$; if they are disjoint, four vertices do.
There are at most $\binom t2$ choices for such a pair of edges. Counting the remaining vertices of the cycle and allowing overcounting, we obtain, for $t\ge4$,
\[
N(C_t,F)\le\binom t2\left(n^{t-3}\sum_{b\in B}d(b)^3+n^{t-4}\sum_{b\in B}d(b)^4\right).
\]
When $t=3$, only the first term occurs. Since $\Delta(\Gamma)\le\delta n$, part (i) gives
\[
\sum_{b\in B}d(b)^3\le c_t\delta n^3,\qquad \sum_{b\in B}d(b)^4\le c_t\delta^2n^4.
\]
Since $\delta\le1/2$, it follows that
\[
N(C_t,F)\le t(t-1)c_t\delta n^t\le\rho n^t.
\]

Applying Lemma~\ref{lem:removal} with the already fixed $\eta$ and $\rho$ (and adding isolated vertices to $F$ if necessary), we obtain $R\subseteq E(F)$ with $|R|\le\eta n^2$ such that $F-R$ is $C_t$-free.

For each $b\in B$, the neighbourhood $N_\Gamma(b)$ spans a clique in $F$. Since $F-R$ is $C_t$-free, the graph induced by $N_\Gamma(b)$ in $F-R$ has at most $\operatorname{ex}(d(b),C_t)$ edges. By our choice of $c_t$, we have $\operatorname{ex}(d(b),C_t)\le 3d(b)^2/8$ whenever $d(b)\ge c_t$. For such vertices,
\[
|R\cap E(F[N_\Gamma(b)])|\ge\binom{d(b)}2-\operatorname{ex}(d(b),C_t)\ge\frac18d(b)^2-\frac12d(b).
\]
For $d(b)<c_t$, we have $d(b)^2\le c_td(b)$. Since $c_t\ge5$, combining these two cases yields, for every $b\in B$,
\[
d(b)^2\le8|R\cap E(F[N_\Gamma(b)])|+c_td(b).
\]

Summing over $b\in B$, we have
\[
\sum_{b\in B}|R\cap E(F[N_\Gamma(b)])|=\sum_{xy\in R}q(x,y).
\]
Combining this with \eqref{eq:codegree} and the Cauchy--Schwarz inequality, we obtain
\[
\begin{aligned}
\sum_{b\in B}d(b)^2
&\le8\sum_{xy\in R}q(x,y)+c_t\sum_{b\in B}d(b)\\
&\le8\sqrt{|R|\sum_{\{x,y\}\subseteq A}q(x,y)^2}+c_te(\Gamma)\\
&\le8\sqrt{c_t\eta}\,n^2+c_te(\Gamma).
\end{aligned}
\]
Interchanging $A$ and $B$ and adding the resulting inequalities, we obtain
\[
\sum_{v\in V(\Gamma)}d(v)^2\le16\sqrt{c_t\eta}\,n^2+2c_te(\Gamma).
\]
Finally, since $2e(\Gamma)\le n\Delta(\Gamma)\le\delta n^2$, the choices of $\eta$ and $\delta$ give
\[
\sum_{v\in V(\Gamma)}d(v)^2\le\left(16\sqrt{c_t\eta}+c_t\delta\right)n^2\le\varepsilon n^2.\qedhere
\]
\end{proof}

\begin{lem}\label{lem:wMantel}
Let $t\ge2$ be an integer, and let $D$ be an oriented graph on $n$ vertices with an arc weight function $w:A(D)\to\{0,1,\ldots,t-1\}$. Suppose that $w(x,y)+w(y,z)\le t-1$ whenever $xyz$ is a directed path in $D$. Then
$\sum_{e\in A(D)}w(e)\le \frac{t-1}{4}n^2$.
\end{lem}

\begin{proof}
For each integer $1\le i\le t-1$, let $E_i=\{e\in A(D):w(e)\ge i\}$. Every arc $e$ belongs to exactly $w(e)$ of these sets, so
\begin{equation}\label{eq:weight-layers}
\sum_{e\in A(D)}w(e)=\sum_{i=1}^{t-1}|E_i|.
\end{equation}

We claim that for each $1\le i\le\lfloor(t-1)/2\rfloor$,
\begin{equation}\label{eq:paired-layers}
|E_i|+|E_{t-i}|\le\frac{n^2}{2}.
\end{equation}
Fix such an $i$. Since $2(t-i)>t-1$, $E_{t-i}$ contains no directed path of length two.
Let $A$ be the set of tails of arcs in $E_{t-i}$ and $B$ the set of their heads. Then $A$ and $B$ are disjoint.
If $uv\in E_{t-i}$, then any arc entering $u$ or leaving $v$ has weight at most $i-1$. Consequently, no arc of $E_i$ enters $A$, and no arc of $E_i$ leaves $B$.
Put $R=V(D)\setminus(A\cup B)$. Thus every arc of $E_i$ goes from $A$ to $B$, from $A$ to $R$, from $R$ to $B$, or has both endpoints in $R$. Since $D$ is oriented, we have
$|E_i|\le |A||B|+|A||R|+|R||B|+\binom{|R|}{2}$.
All arcs of $E_{t-i}$ go from $A$ to $B$, so $|E_{t-i}|\le |A||B|$. Since $|A|+|B|+|R|=n$, we obtain
\[
|E_i|+|E_{t-i}|
\le 2|A||B|+(|A|+|B|)|R|+\frac{|R|^2}{2}
=\frac{n^2-(|A|-|B|)^2}{2}
\le\frac{n^2}{2}.
\]

If $t$ is odd, the $t-1$ levels can be paired as $(E_i,E_{t-i})$ for $1\le i\le (t-1)/2$.
Summing \eqref{eq:paired-layers} over these pairs and using \eqref{eq:weight-layers}, we obtain
$\sum_{e\in A(D)}w(e)\le\frac{t-1}{4}n^2$.

If $t$ is even, the levels other than $E_{t/2}$ can be paired as $(E_i,E_{t-i})$ for $1\le i\le (t-2)/2$. Since two consecutive arcs in $E_{t/2}$ would have total weight at least $t$, $E_{t/2}$ contains no directed path of length two.
Let $P$ and $Q$ be its sets of tails and heads. As above, $P\cap Q=\emptyset$ and every arc of $E_{t/2}$ goes from $P$ to $Q$. Therefore $|E_{t/2}|\le |P||Q| \le\frac{n^2}{4}$.
Together with \eqref{eq:weight-layers} and \eqref{eq:paired-layers}, this yields $\sum_{e\in A(D)}w(e)\le\frac{t-1}{4}n^2$.
\end{proof}

We have now proved the necessary lemmas and are ready to establish the main theorem.

\begin{proof}[\bf Proof of Theorem~\ref{thm:main}]
For the lower bound, let $V_1,V_2,V_3$ be three disjoint independent sets with $|V_1|=t-1$ and $|V_2|,|V_3|$ as equal as possible, subject to $|V_1|+|V_2|+|V_3|=n$. Add all arcs from $V_1$ to $V_2$, from $V_2$ to $V_3$, and from $V_3$ to $V_1$. The resulting oriented graph $G$ is $\Cdir_{3t}$-free, since every directed cycle visits the three parts equally often, whereas $|V_1|=t-1$. Moreover, every directed triangle has one vertex in each part, so
\[
N(\Cdir_3,G)=|V_1||V_2||V_3|=\frac{t-1}{4}n^2-O_t(n).
\]

It remains to establish the upper bound. Let $G$ be an $n$-vertex $\Cdir_{3t}$-free oriented graph, where $t\ge 3$ is an integer. For an arc $uv\in A(G)$, let $\tau_G(uv)$ denote the number of directed triangles in $G$ containing $uv$. An arc $uv$ is called \emph{heavy} if $\tau_G(uv)\ge 3t-4$ and \emph{light} otherwise. Let $H$ and $L$ be the spanning subgraphs of $G$ consisting of its heavy and light arcs, respectively.
Let $T_j(G)$ denote the number of directed triangles in $G$ containing exactly $j$ heavy arcs, where $j\in\{0,1,2,3\}$. Then
\[
\sum_{uv\in A(G)}\tau_G(uv)=3N(\Cdir_3,G)
\text{~~and~~}
N(\Cdir_3,G)=\sum_{j=0}^3T_j(G).
\]

We first bound the number of arcs in $H$.

\begin{claim}\label{lem:alternating}
For any disjoint $X,Y\subseteq V(G)$, the underlying graph of the subgraph of $H$ consisting of all arcs from $X$ to $Y$ is $C_{2t}$-free. Consequently, $e(H)=O_t(n^{1+1/t})$.
\end{claim}

\begin{proof}[\bf Proof of Claim~\ref{lem:alternating}]
Suppose to the contrary that this bipartite graph contains a cycle
$a_1b_1a_2b_2\cdots a_t\\b_ta_1$, where $a_i\in X$, $b_i\in Y$, and indices are taken modulo $t$. Since all its edges correspond to heavy arcs from $X$ to $Y$, the arcs $a_ib_i$ and $a_{i+1}b_i$ are heavy for every $i\in[t]$.
For each $i$, a directed triangle in $G$ containing the heavy arc $a_{i+1}b_i$ has the form $a_{i+1}b_ic_ia_{i+1}$. The vertex $c_i$ cannot be $a_{i+1}$ or $b_i$, since these are already vertices of the triangle. It cannot be $a_i$ or $b_{i+1}$ either, since this would require the arc $b_ia_i$ or $b_{i+1}a_{i+1}$, respectively.
Thus at most $2t-4$ vertices on the cycle can serve as $c_i$.
Since each $a_{i+1}b_i$ is heavy, it belongs to at least $3t-4$ directed triangles, so there are at least $t$ choices for $c_i$ outside the cycle. We may greedily choose pairwise distinct vertices $c_1,\ldots,c_t$ such that $a_{i+1}b_ic_ia_{i+1}$ is a directed triangle for every $i$. Then
$a_1b_1c_1a_2b_2c_2\cdots a_tb_tc_ta_1$ is a directed cycle of length $3t$ in $G$, contradicting the assumption that $G$ is $\Cdir_{3t}$-free.

To bound $e(H)$, randomly partition $V(G)$ into two sets $X$ and $Y$ by assigning each vertex independently and uniformly to one of them.
Let $F$ be the undirected bipartite graph with bipartition $(X,Y)$, where $xy\in E(F)$ if $xy\in A(H)$ for $x\in X$ and $y\in Y$. Since $F$ is $C_{2t}$-free, Lemma~\ref{lem:BS} gives $e(F)\le 100t\,n^{1+1/t}$. Each arc $xy\in A(H)$ appears as an edge of $F$ precisely when $x\in X$ and $y\in Y$, which occurs with probability $1/4$. Thus, by linearity of expectation,
\[
100t\,n^{1+1/t}\ge \mathbb{E}[e(F)]
=\sum_{xy\in A(H)}\mathbb{P}(x\in X,\ y\in Y)
=\frac14 e(H).
\]
Hence $e(H)=O_t(n^{1+1/t})$.
\end{proof}

We next show that directed triangles containing at most one heavy arc contribute only $o(n^2)$.

\begin{claim}\label{lem:few0}
$T_0(G)=o(n^2)$.
\end{claim}

\begin{proof}[\bf Proof of Claim~\ref{lem:few0}]
Every arc of $L$ belongs to at most $3t-5$ directed triangles. Then,
\[
3T_0(G)=\sum_{uv\in A(L)}\tau_L(uv)
\le (3t-5)e(L)=O_t(n^2)=o(n^3).
\]
Fix $\eta>0$, and let $\xi>0$ be given by Lemma~\ref{lem:dirremoval}. For sufficiently large $n$, we have $T_0(G)\le\xi n^3$, so there is a set $R\subseteq A(L)$ with $|R|\le\eta n^2$ such that $L-R$ is $\Cdir_3$-free.
Every directed triangle in $L$ contains at least one arc of $R$, and each arc of $R$ belongs to at most $3t-5$ directed triangles in $L$. Hence $T_0(G)\le(3t-5)|R|\le(3t-5)\eta n^2$. Since $\eta>0$ is arbitrary, $T_0(G)=o(n^2)$.
\end{proof}

\begin{claim}\label{lem:few1}
$T_1(G)=O_t\!\left(n^{\frac32+\frac1t+\frac1{2t^2}}\right)=o(n^2).$
\end{claim}
\begin{proof}[\bf Proof of Claim~\ref{lem:few1}]
Fix a partition $V(G)=X\cup Y\cup Z$. Let $Q$ be the undirected bipartite graph with bipartition $(U,Z)$, where $U=\{(a,b):a\in X,\ b\in Y,\ ab\in A(H)\}$. Join $(a,b)\in U$ to $c\in Z$ if $abc$ is a directed triangle in $G$ and both $bc$ and $ca$ are light. Write $m=|U|$ and $r=|Z|$. By Claim~\ref{lem:alternating}, we have $m\le e(H)=O_t(n^{1+1/t})$.
We aim to show that $e(Q)=O_t\!\left(n^{\frac32+\frac1t+\frac1{2t^2}}\right)$. This is trivial if $e(Q)=0$, so we may assume $e(Q)>0$. In particular, $m,r>0$.

Let $A(Q)$ be the adjacency matrix of $Q$, and let $\lambda_1,\ldots,\lambda_{m+r}$ be its eigenvalues. Since $A(Q)$ is symmetric and nonnegative, the Perron--Frobenius theorem implies that $\rho(Q)$ is the largest eigenvalue of $A(Q)$. Since homomorphisms from $C_{2t}$ to $Q$ correspond to closed walks of length $2t$, we have
\[
\hom(C_{2t},Q) =\operatorname{tr}(A(Q)^{2t}) =\sum_{i=1}^{m+r}\lambda_i^{2t} \ge 2\rho(Q)^{2t},
\]
where the inequality follows from the symmetry of the spectrum of the bipartite graph $Q$.
Let $\boldsymbol{1}_m$ and $\boldsymbol{1}_r$ denote the all-ones vectors of dimensions $m$ and $r$, respectively. Consider the vector
\[
\mathbf{x}= \left( \sqrt{r}\,\boldsymbol{1}_m^{\mathsf T}, \sqrt{m}\,\boldsymbol{1}_r^{\mathsf T} \right)^{\mathsf T},
\]
where the vertices of $Q$ are ordered with those in $U$ preceding those in $Z$. By Lemma~\ref{lem:Rayleigh},
\[
\rho(Q)\ge\frac{\mathbf{x}^{\mathsf T}A(Q)\mathbf{x}} {\mathbf{x}^{\mathsf T}\mathbf{x}}
=\frac{2\sum_{uv\in E(Q)}x_ux_v}{2mr}
=\frac{e(Q)}{\sqrt{mr}}.
\]
Consequently,
\begin{equation}\label{eq:hom-lower}
\hom(C_{2t},Q)^{1/(2t)} \ge \frac{2^{1/(2t)}e(Q)}{\sqrt{mr}}.
\end{equation}

Define a symmetric relation $\sim$ on $V(Q)$ as follows: for $(a,b),(a',b')\in U$, $(a,b)\sim(a',b')$ if and only if $a=a'$ or $b=b'$; for $c,c'\in Z$, $c\sim c'$ if and only if $c=c'$; and no vertex of $U$ is related to a vertex of $Z$.

We verify that $|\{w\in N_Q(v):u\sim w\}|\le 2(3t-5)$ for all $u,v\in V(Q)$. If $u$ and $v$ lie in the same part of $Q$, the set is empty. If $u\in Z$ and $v\in U$, the set contains at most one vertex.
Now suppose $u=(a,b)\in U$ and $v=c\in Z$. Every neighbour $w=(a',b')$ of $c$ satisfying $u\sim w$ has $a'=a$ or $b'=b$. If $a'=a$, then $ab'c$ is a directed triangle in $G$, and $ca$ is light. Hence there are at most $\tau_G(ca)\le 3t-5$ choices for $b'$. Similarly, if $b'=b$, then $a'bc$ is a directed triangle in $G$, and $bc$ is light, giving at most $\tau_G(bc)\le 3t-5$ choices for $a'$. Therefore,
$|\{w\in N_Q(c):(a,b)\sim w\}|\le 2(3t-5)$.

We also have $\Delta(Q)\le (3t-5)n$. Indeed, each vertex of $U$ has at most $|Z|\le n$ neighbours. Now fix $c\in Z$. For each $a\in X$, consider the vertices $b\in Y$ such that $(a,b)$ is adjacent to $c$ in $Q$. If there are any such vertices, then $ca$ is a light arc of $G$, and each such $b$ completes $ca$ to a directed triangle. Hence there are at most $\tau_G(ca)\le 3t-5$ choices of $b$. Summing over $a\in X$, we obtain
$d_Q(c)\le (3t-5)|X|\le (3t-5)n.$

We claim that for every homomorphism $f:C_{2t}\to Q$, there exist distinct vertices $x,y\in V(C_{2t})$ such that $f(x)\sim f(y)$.
Suppose, for a contradiction, that there is a homomorphism $f:C_{2t}\to Q$ satisfying $f(x)\not\sim f(y)$ for all distinct $x,y\in V(C_{2t})$.
Since $Q$ is bipartite, the homomorphism $f$ gives a closed walk
$(a_1,b_1),c_1,(a_2,b_2),c_2,\ldots,(a_t,b_t),c_t$, where $(a_i,b_i)\in U$ and $c_i\in Z$.
By our assumption and the definition of $\sim$, the vertices $a_1,\ldots,a_t$ are pairwise distinct, as are $b_1,\ldots,b_t$ and $c_1,\ldots,c_t$. Since $X,Y,Z$ are disjoint, all these $3t$ vertices are distinct.
By the definition of $Q$, $a_1b_1c_1a_2b_2c_2\cdots a_tb_tc_ta_1$ is a directed cycle of length $3t$ in $G$, a contradiction.

Since every homomorphism from $C_{2t}$ to $Q$ satisfies the relation condition in Lemma~\ref{lem:janzer}, applying the lemma with $k=t$ and $s=2(3t-5)$ gives
\[
\hom(C_{2t},Q) \le 32t^{3/2}\sqrt{2(3t-5)\Delta(Q)} (m+r)^{1/(2t)} \hom(C_{2t},Q)^{1-1/(2t)}.
\]
Since $e(Q)>0$, we have $\hom(C_{2t},Q)>0$. Thus,
\[
\hom(C_{2t},Q)^{1/(2t)} \le 32t^{3/2}\sqrt{2(3t-5)\Delta(Q)} (m+r)^{1/(2t)}.
\]
Combining this with \eqref{eq:hom-lower} and using $m=O_t(n^{1+1/t})$, $r\le n$, and $\Delta(Q)=O_t(n)$,
we obtain
\[
e(Q) \le \frac{32t^{3/2}\sqrt{2(3t-5)}}{2^{1/(2t)}} \sqrt{mr\Delta(Q)}(m+r)^{1/(2t)}
=O_t\!\left(n^{\frac32+\frac1t+\frac1{2t^2}}\right).
\]

Finally, randomly partition $V(G)$ into three sets $X,Y,Z$ by assigning each vertex independently and uniformly to one of them.
Every directed triangle with exactly one heavy arc has a unique cyclic ordering $(a,b,c)$ in which $ab$ is heavy. By the definition of $Q$, such a triangle contributes exactly one edge to $Q$ if and only if $a\in X$, $b\in Y$, and $c\in Z$.
Thus $e(Q)$ counts precisely the triangles satisfying these conditions. By linearity of expectation,
\[
\mathbb{E}[e(Q)]
=\sum_{(a,b,c)}\mathbb{P}(a\in X,\ b\in Y,\ c\in Z)
=\frac{T_1(G)}{27},
\]
where the sum is over all directed triangles with exactly one heavy arc, written in the above cyclic ordering.
Since the bound on $e(Q)$ holds for every partition, we conclude that
\[
T_1(G)=O_t\!\left(n^{\frac32+\frac1t+\frac1{2t^2}}\right)=o(n^2).\qedhere
\]
\end{proof}

We next show that, after removing a bounded number of vertices, the number of remaining directed triangles can be made arbitrarily small compared with $n^2$.

\begin{claim}\label{lem:cover}
For every $\varepsilon>0$, there exist positive integers $K=K(\varepsilon,t)$ and $n_0=n_0(\varepsilon,t)$ such that every $\Cdir_{3t}$-free oriented graph $G$ on $n\ge n_0$ vertices contains a set $S\subseteq V(G)$ satisfying $|S|\le K$ and $N(\Cdir_3,G-S)\le \varepsilon n^2$.
\end{claim}

\begin{proof}[\bf Proof of Claim~\ref{lem:cover}]
Fix $\varepsilon>0$. Let $\delta=\delta(t,\varepsilon/8)$ be given by Lemma~\ref{lem:dsq}~(ii), and set $K=\left\lceil\frac{16c_t}{\delta^2}\right\rceil$. By Claims~\ref{lem:few0} and~\ref{lem:few1}, there exists a positive integer $n_0=n_0(\varepsilon,t)$ such that $T_0(G)+T_1(G)\le\frac{\varepsilon}{2}n^2$ for every $\Cdir_{3t}$-free oriented graph $G$ on $n\ge n_0$ vertices.

Now let $G$ be such a graph. Randomly partition $V(G)$ into two sets $X,Y$ by assigning each vertex independently and uniformly to one of them.
Let $\Gamma$ be the undirected bipartite graph with bipartition $(X,Y)$, where $xy\in E(\Gamma)$ if $x\in X$, $y\in Y$, and $xy\in A(H)$. By Claim~\ref{lem:alternating}, $\Gamma$ is $C_{2t}$-free. Thus, Lemma~\ref{lem:dsq}~(i) gives $\sum_{v\in V(G)}d_\Gamma(v)^2\le c_tn^2$ for every choice of $(X,Y)$.

Let $v\in V(G)$. For each arc $e$ of $H$ incident with $v$, let $I_e$ be the indicator of the event that the underlying edge of $e$ belongs to $E(\Gamma)$. Then
\[
d_\Gamma(v)^2
=\Bigl(\sum_{e\ni v}I_e\Bigr)^2
=\sum_{e\ni v}I_e +\sum_{\substack{e,f\ni v\\e\ne f}}I_eI_f.
\]
For each arc $e$ incident with $v$, we have $\mathbb E[I_e]=1/4$. For any two distinct arcs $e,f$ incident with $v$,
\[
\mathbb E[I_eI_f]=
\begin{cases}
\mathbb P(v\in X,\ u,w\in Y)=\frac18,
& \text{if $e=vu$ and $f=vw$},\\
\mathbb P(v\in Y,\ u,w\in X)=\frac18,
& \text{if $e=uv$ and $f=wv$},\\
\mathbb P(v\in X\cap Y)=0,
& \text{if one is incoming and the other outgoing}.
\end{cases}
\]
Taking expectations, we obtain
\[
\begin{aligned}
\mathbb E[d_\Gamma(v)^2]
&=\frac14\bigl(d_H^+(v)+d_H^-(v)\bigr)+\frac18\bigl[d_H^+(v)(d_H^+(v)-1)+d_H^-(v)(d_H^-(v)-1)\bigr]\\
&=\frac18\bigl(d_H^+(v)^2+d_H^-(v)^2\bigr)+\frac18\bigl(d_H^+(v)+d_H^-(v)\bigr).
\end{aligned}
\]
Summing over all $v\in V(G)$ and using linearity of expectation and
$\sum_{v\in V(G)}(d_H^+(v)+d_H^-(v))=2e(H)$, we obtain
\[
\frac18\sum_{v\in V(G)} \bigl(d_H^+(v)^2+d_H^-(v)^2\bigr)+\frac14e(H)
=\sum_{v\in V(G)}\mathbb E[d_\Gamma(v)^2]
=\mathbb E\Bigl[\sum_{v\in V(G)}d_\Gamma(v)^2\Bigr]
\le c_tn^2.
\]
Thus,
\begin{equation}\label{eq:heavy-square}
\sum_{v\in V(G)}\bigl(d_H^+(v)^2+d_H^-(v)^2\bigr) \le 8c_tn^2.
\end{equation}

Let $S=\{v\in V(G):d_H^+(v)+d_H^-(v)>\delta n\}$ and $W=V(G)\setminus S$.
Since $d_H^+(v)^2+d_H^-(v)^2\ge\frac12(d_H^+(v)+d_H^-(v))^2$, it follows from \eqref{eq:heavy-square} that
\[
\frac12|S|\delta^2n^2
\le\sum_{v\in S}\bigl(d_H^+(v)^2+d_H^-(v)^2\bigr)
\le 8c_tn^2.
\]
Thus, $|S|\le\frac{16c_t}{\delta^2}\le K$. 

Apply the same random bipartition procedure to $H[W]$, and let $\Gamma_W$ be the resulting undirected bipartite graph. By Claim~\ref{lem:alternating}, $\Gamma_W$ is $C_{2t}$-free. Moreover, $\Delta(\Gamma_W)\le\delta n$ by the definition of $W$.
Adding $n-|W|$ isolated vertices to $\Gamma_W$ preserves its $C_{2t}$-freeness. Hence, Lemma~\ref{lem:dsq}~(ii) gives $\sum_{v\in W}d_{\Gamma_W}(v)^2 \le\frac{\varepsilon}{8}n^2$. Taking expectations and repeating the preceding calculation for $H[W]$, we obtain
\[
\frac18\sum_{v\in W}
\bigl(d_{H[W]}^+(v)^2+d_{H[W]}^-(v)^2\bigr)+\frac14e(H[W])
=\mathbb E\left[\sum_{v\in W}d_{\Gamma_W}(v)^2\right]
\le\frac{\varepsilon}{8}n^2.
\]
Thus,
\begin{equation}\label{eq:remaining-heavy-square}
\sum_{v\in W}
\bigl(d_{H[W]}^+(v)^2+d_{H[W]}^-(v)^2\bigr)
\le\varepsilon n^2.
\end{equation}

Every directed triangle in $G[W]$ containing at least two heavy arcs contains a directed path of length two in $H[W]$. Since each such path determines at most one directed triangle, the number of these triangles is at most
\[
\sum_{v\in W}d_{H[W]}^-(v)d_{H[W]}^+(v)
\le\frac12\sum_{v\in W}\bigl(d_{H[W]}^+(v)^2+d_{H[W]}^-(v)^2\bigr)
\le\frac{\varepsilon}{2}n^2,
\]
where the last inequality follows from \eqref{eq:remaining-heavy-square}. The number of directed triangles in $G[W]$ containing at most one heavy arc is at most $T_0(G)+T_1(G)\le\frac{\varepsilon}{2}n^2$. Therefore,
\[
N(\Cdir_3,G-S)=N(\Cdir_3,G[W])\le\varepsilon n^2.\qedhere
\]
\end{proof}

We now bound the number of directed triangles meeting a given vertex set. For an oriented graph $G$ and a set $S\subseteq V(G)$, let $N_S(\Cdir_3,G)$ denote the number of directed triangles in $G$ containing at least one vertex of $S$.

\begin{claim}\label{lem:meet}
Let $G$ be a $\Cdir_{3t}$-free oriented graph on $n$ vertices, and let $S\subseteq V(G)$ with $|S|=k$. Then
\[
N_S(\Cdir_3,G) \le \frac{t-1}{4}(n-k)^2 +2(3t-1)k\,9^k n +\binom{k}{2}(n-k) +\binom{k}{3}.
\]
\end{claim}

\begin{proof}[\bf Proof of Claim~\ref{lem:meet}]
Put $W=V(G)\setminus S$. We say that two vertices $u,v\in W$ are equivalent if
$N_S^-(u)=N_S^-(v)$ and $N_S^+(u)=N_S^+(v)$.
Partition $W$ into equivalence classes, which we call types. Since $G$ is oriented, we have
$N_S^-(v)\cap N_S^+(v)=\emptyset$ for every $v\in W$. The number of ordered pairs $(S_1,S_2)$ of disjoint subsets of $S$ is $\sum_{i=0}^{k}\binom{k}{i}2^{k-i}=3^k$. Hence $W$ has at most $3^k$ types.

For each arc $uv\in A(G[W])$, assign the weight $w(u,v)=|N_S^-(u)\cap N_S^+(v)|$. For every $s\in N_S^-(u)\cap N_S^+(v)$, $suv$ is a directed triangle. 
Then the number of directed triangles containing exactly one vertex of $S$ is
$\sum_{uv\in A(G[W])}w(u,v)$.
The number of directed triangles containing at least two vertices of $S$ is at most
$\binom{k}{2}(n-k)+\binom{k}{3}$.

For any two types $X,Y$, define $w(X,Y)=|N_S^-(x)\cap N_S^+(y)|$, where $x\in X$ and $y\in Y$. All arcs from $X$ to $Y$ have weight $w(X,Y)$. We call $(X,Y)$ an ordered type pair. There are at most $9^k$ such pairs, and every arc has weight at most $k$.
For each ordered type pair $(X,Y)$, discard all arcs from $X$ to $Y$ if they contain no matching of size $3t$, and retain them otherwise. Let $F$ be the weighted oriented graph on $W$ consisting of the retained arcs.

We first show that the total weight of the discarded arcs is at most
\begin{equation}\label{eq:discarded-weight}
2(3t-1)k\,9^k n.
\end{equation}
Indeed, fix a discarded ordered type pair $(X,Y)$, and let $U$ be the set of endpoints of a maximum matching $M$ among the arcs from $X$ to $Y$.
Then $|U|\le 2(3t-1)$. By the maximality of $M$, every arc from $X$ to $Y$ has an endpoint in $U$. Hence the total weight of the discarded arcs from this pair is at most
$w(X,Y)|U|n\le 2(3t-1)kn$. Summing over at most $9^k$ ordered type pairs proves \eqref{eq:discarded-weight}.

We next show that
\begin{equation}\label{eq:arc-weight}
0\le w(e)\le t-1 \text{~~for~every~~} e\in A(F).
\end{equation}
Suppose, to the contrary, that $F$ contains an arc $xy$ with $w(x,y)\ge t$. Let $X$ and $Y$ be the types containing $x$ and $y$, respectively. Since $(X,Y)$ is retained in $F$, there exist $t$ pairwise vertex-disjoint arcs $x_1y_1,\ldots,x_ty_t$ in $F$, with $x_i\in X$ and $y_i\in Y$ for every $i\in[t]$. Since $w(X,Y)\ge t$, we may choose distinct vertices $s_1,\ldots,s_t\in N_S^-(x)\cap N_S^+(y)$. By the definition of types, $s_j\in N_S^-(x_i)\cap N_S^+(y_i)$ for all $i,j\in\{1,\ldots,t\}$. Thus
$s_1x_1y_1s_2x_2y_2\cdots s_tx_ty_ts_1$ is a directed cycle of length $3t$, contradicting the assumption that $G$ is $\Cdir_{3t}$-free.

We now prove that every directed path $xyz$ in $F$ satisfies
\begin{equation}\label{eq:two-path-weight}
w(x,y)+w(y,z)\le t-1.
\end{equation}
Suppose, to the contrary, that $w(x,y)+w(y,z)\ge t$.
Let $X,Y,Z$ be the types containing $x,y,z$, respectively. By \eqref{eq:arc-weight},
$w(X,Y)\le t-1$ and $w(Y,Z)\le t-1$.
Since $w(X,Y)+w(Y,Z)=w(x,y)+w(y,z)\ge t$, both weights are positive. Hence we can choose positive integers $\alpha,\beta$ such that
\[
\alpha+\beta=t, ~~\alpha\le w(X,Y), ~~\beta\le w(Y,Z).
\]
Choose distinct vertices $l_1,\ldots,l_\alpha\in N_S^-(x)\cap N_S^+(y)$ and distinct vertices $r_1,\ldots,r_\beta\in N_S^-(y)\cap N_S^+(z)$. These $t$ vertices are distinct since $N_S^+(y)\cap N_S^-(y)=\emptyset$.
Since both ordered type pairs $(X,Y)$ and $(Y,Z)$ are retained in $F$, the arcs corresponding to each pair contain a matching of size at least $3t$.
We may therefore greedily choose $\alpha-1$ arcs $x_iy_i$ from $X$ to $Y$ and $\beta-1$ arcs $y'_jz'_j$ from $Y$ to $Z$, all with pairwise distinct endpoints outside $\{x,y,z\}$.
Since $w(X,Y),w(Y,Z)>0$ and $w(T,T)=0$ for every type $T$, we have $X\ne Y$ and $Y\ne Z$. Thus exactly $1+(\alpha-1)+(\beta-1)=t-1$ of the previously chosen vertices of $W$ lie in $Y$. As $|Y|\ge 3t$, we may choose a further vertex $y^*\in Y$ distinct from all previously chosen vertices in $Y$.
Now consider the directed cycle
\[
l_1x_1y_1l_2\cdots x_{\alpha-1}y_{\alpha-1}l_\alpha xyzr_1y'_1z'_1r_2\cdots y'_{\beta-1}z'_{\beta-1}r_\beta y^*l_1.
\]
The corresponding sequences of auxiliary arcs are omitted when $\alpha=1$ or $\beta=1$. The cycle has $t$ distinct vertices in $S$ and $2(\alpha-1)+3+2(\beta-1)+1=2t$ distinct vertices in $W$, giving a directed cycle of length $3t$, a contradiction. This proves \eqref{eq:two-path-weight}.

By \eqref{eq:arc-weight} and \eqref{eq:two-path-weight}, Lemma~\ref{lem:wMantel} applied to $F$ gives
\[
\sum_{e\in A(F)}w(e) \le \frac{t-1}{4}|W|^2 =\frac{t-1}{4}(n-k)^2.
\]
Adding the discarded weight from \eqref{eq:discarded-weight} and the number of directed triangles containing at least two vertices of $S$, we obtain
\[
N_S(\Cdir_3,G) \le \frac{t-1}{4}(n-k)^2 +2(3t-1)k\, 9^k n +\binom{k}{2}(n-k)+\binom{k}{3}.\qedhere
\]
\end{proof}

Fix $\varepsilon>0$. By Claim~\ref{lem:cover}, there exist constants $K=K(\varepsilon,t)$ and $n_0=n_0(\varepsilon,t)$ such that, for every $n\ge n_0$ and every $\Cdir_{3t}$-free oriented graph $G$ on $n$ vertices, there is a set $S\subseteq V(G)$ with $|S|\le K$ and
$N(\Cdir_3,G-S)\le \varepsilon n^2$.
By Claim~\ref{lem:meet} and $|S|\le K$, we have
\[
N_S(\Cdir_3,G) \le \frac{t-1}{4}n^2 +2(3t-1)K\,9^K n +\binom{K}{2}n+\binom{K}{3} =\frac{t-1}{4}n^2+O_{t,K}(n).
\]
Consequently,
\[
N(\Cdir_3,G) =N(\Cdir_3,G-S)+N_S(\Cdir_3,G) \le \frac{t-1}{4}n^2+\varepsilon n^2+O_{t,K}(n).
\]
Since $K$ depends only on $t$ and $\varepsilon$, and $\varepsilon>0$ is arbitrary, it follows that
\[
\operatorname{ex}(n,\Cdir_3,\Cdir_{3t}) \le \frac{t-1}{4}n^2+o(n^2).\qedhere
\]
\end{proof}

\section{\normalsize Concluding Remarks}\label{s4}

Grzesik et al.~\cite{GJ} determined the order of magnitude of $\ex(n,\Cdir_k,\Cdir_\ell)$ for all fixed distinct integers $k,\ell\ge3$, and obtained asymptotically sharp results in several cases where $k\nmid\ell$. However, when $k\mid\ell$, determining the leading term remains open in general. Our main result, together with the case $t=2$ established in~\cite{GJ}, determines
\[
\ex(n,\Cdir_3,\Cdir_{3t})
=\frac{t-1}{4}n^2+o(n^2)
\]
for every fixed integer $t\ge2$, thereby settling the case $k=3$. An asymptotically optimal construction is a blow-up of $\Cdir_3$ with one class of size $t-1$ and two classes whose sizes differ by at most one.

More generally, for fixed integers $k\ge3$ and $t\ge2$, a blow-up of $\Cdir_k$ with one class of size $t-1$ and the remaining $k-1$ classes as equal as possible gives
\[
\ex(n,\Cdir_k,\Cdir_{kt}) \ge \frac{t-1}{(k-1)^{k-1}}\,n^{k-1} - O_{k,t}(n^{k-2}).
\]
This elementary construction provides a general lower bound. For $k\ge4$, the main challenge is to determine the sharp leading coefficient and identify the corresponding asymptotically extremal configurations.

The case $k\mid\ell$ requires finer counting arguments than the case $k\nmid\ell$, since the extremal number is of order $n^{k-1}$ rather than $n^k$. In general, an arc may belong to $O_k(n^{k-2})$ directed $k$-cycles, so deleting $o(n^2)$ arcs gives only an $o(n^k)$ bound on the number of cycles destroyed, which is insufficient to determine the leading term at the $n^{k-1}$ scale. In our proof for $k=3$, we overcome this difficulty in part by distinguishing heavy and light arcs and using the bounded number of directed triangles containing each light arc. Extending such refined counting and structural arguments to longer directed cycles remains a challenge. This motivates the following problem.

\begin{pb}
For fixed integers $k\ge4$ and $t\ge2$, determine the asymptotic value of $\operatorname{ex}(n,\Cdir_k,\Cdir_{kt})$.
\end{pb}

\end{document}